\documentclass[a4paper,11pt]{article}
\pdfoutput=1 

\usepackage{jheppub} 

\usepackage[T1]{fontenc} 

\title{\boldmath The Enriquez connection for higher genus polylogarithms}

\author{Takashi Ichikawa}
\affiliation{Department of Mathematics, Faculty of Science and Engineering, 
Saga University, 
\\ 
Saga 840-8502, Japan}

\emailAdd{ichikawn@cc.saga-u.ac.jp}

\abstract{
We study the variation of the Enriquez connection for higher genus polylogarithms 
under degenerations of Riemann surfaces with marked points, 
and show that this connection becomes the connection constructed by the author 
for degenerating families of pointed Riemann surfaces. 
Therefore, 
we have an important application that the higher genus polylogarithms 
derived from the Enriquez connection can be described explicitly as power series 
in deformation parameters and their logarithms associated with the families 
whose coefficients are expressed by multiple zeta values. 
\vspace{4ex}

\noindent
{\sc Keywords:} Differential and Algebraic Geometry, Scattering Amplitudes 
\vspace{4ex}

\noindent
{\sc ArXiv ePrint: 2510.00486} 
}

\begin{document} 
\maketitle
\flushbottom

\section{Introduction}

{\it Polylogarithm functions}, or {\it polylogarithms} for short, 
describe monodromies of nilpotent connections defined on families of algebraic varieties. 
Especially, polylogarithms on families of pointed Riemann surfaces become central objects 
for studying Feynman integrals in quantum field theory and motives in arithmetic geometry. 
For the references, see \cite{GiEG} for the genus $0$ case, 
and \cite{BMW, BB, BDDT, BDDPT1, BDDPT2, BDDPT3, BK, MS, W} for the genus $1$ case 
in which the Knizhnik-Zamolodchikov (KZ) connections on pointed Riemann spheres 
and their elliptic extensions known as the elliptic Knizhnik-Zamolodchikov-Bernard (KZB) connections 
(cf. \cite{BrowL, CEE, E2, H}) played important roles. 
Recently, 
D'Hoker, Enriquez, Hidding, Schlotterer and Zerbini constructed higher genus polylogarithms 
as monodromies of a single-valued, but non-meromorphic connection 
on each Riemann surface (cf. \cite{DHS}), 
and related this connection with the meromorphic multi-valued connection 
given by Enriquez \cite{E1} (cf. \cite{DESZ}). 
Furthermore, 
D'Hoker and Schlotterer described the Enriquez connection explicitly 
in terms of Abelian differentials (cf. \cite{DS}). 
See \cite{EZ1, EZ2, EZ3, BBILZ} for related results. 

In this paper, 
we study the variation of the Enriquez connection under degenerations of Riemann surfaces 
with marked points using results of \cite{DS, I1, I2}, 
and show that this connection becomes the connection constructed in \cite{I3} 
for degenerating families of pointed Riemann surfaces. 
Therefore, 
we have the following applications: 
\begin{itemize}

\item 
The connection given in \cite{I3} which is locally defined 
for degenerating families of pointed Riemann surfaces 
can be globally extended to all families of pointed Riemann surfaces. 

\item 
The higher genus polylogarithms derived from the Enriquez connection can be described explicitly 
as power series in deformation parameters and their logarithms associated with the families 
whose coefficients are expressed by multiple zeta values. 
This description is a higher genus extension of results by Brown \cite{Brow}, 
Banks-Panzer-Pym \cite{BPP}, Enriquez \cite{E3} in the genus $\leq 1$ case. 

\item 
Based on results of Enriquez and Zerbini \cite{EZ2, EZ3}, 
our higher genus polylogarithms with their explicit formulas are expected to be useful 
in calculating hyperlogarithms on families of pointed Riemann surfaces 
which were studied in the genus $1$ case by Br\"{o}del and others 
\cite{BDDT, BDDPT1, BDDPT2, BDDPT3, BK} in order to calculate Feynman integrals. 

\end{itemize}

\section{The Enriquez kernels and connection}

We consider a compact Riemann surface $R$ of any genus $g$, 
and denote its universal covering space by $\tilde{R}$, 
the associated projection by $\pi : \tilde{R} \rightarrow R$. 
Let $\{ A_{i}, B_{i} \}_{1 \leq i \leq g}$ be a basis of the first homology cycles which is symplectic, 
namely, their intersections satisfy that $A_{i} \cdot A_{j} = B_{i} \cdot B_{j} = 0$ 
and $A_{i} \cdot B_{j}$ is the Kronecker delta $\delta_{ij}$. 
Then there exists a unique basis $\{ \omega_{i} \}_{1 \leq i \leq g}$ 
of the space of holomorphic Abelian differentials on $R$ which is normalized for $\{ A_{i} \}$ 
in the sense that $\oint_{A_{i}} \omega_{j} = \delta_{ij}$, 
and hence one can define the period matrix $\Omega = (\Omega_{ij})_{i,j}$ 
of the surface $R$ with $\{ A_{i}, B_{i} \}$ by 
$$
\oint_{B_{i}} \omega_{j} = \Omega_{ij}. 
$$ 
Choosing the cycles $A_{i}$ and $B_{i}$ so that they share a common base point $q$, 
we regard these cycles as generators of the fundamental group $\pi_{1}(R, q)$ of $R$ 
satisfying 
$$
[A_{1}, B_{1}]  [A_{2}, B_{2}] \cdots [A_{g}, B_{g}] = 1, 
$$
where $[A_{i}, B_{i}]$ denotes the commutator $A_{i} B_{i} A_{i}^{-1} B_{i}^{-1}$. 
Then one can take a fundamental domain $D$ in $\tilde{R}$ for the action of $\pi_{1}(R, q)$ 
which is obtained by cutting $R$ along these cycles. 

Following \cite{DESZ, DS}, 
we review results on the {\it Enriquez connection} which was given by Enriquez \cite{E1} 
as a meromorphic connection $d - {\mathcal K}_{R}$ on $\tilde{R}$ with only simple poles 
which is valued in the freely completed Lie algebra ${\mathfrak g}$ with generators 
$a_{i}, b_{i}$ corresponding to $A_{i}, B_{i}$ respectively. 
We denote its coefficient functions by  
$$
g^{i_{1} \cdots i_{r}}_{j}(x, y) = (g_{R})^{i_{1} \cdots i_{r}}_{j}(x, y) 
\quad (r \geq 0, \ i_{1},..., i_{r} \in \{ 1,..., g \})
$$
which are related to those used in \cite{E1} by 
$g^{i_{1} \cdots i_{r}}_{j}(x, y) = (-2 \pi \sqrt{-1})^{r} \omega^{i_{1} \cdots i_{r}}_{j}(x, y)$. 
These coefficient functions, also referred to as {\it Enriquez kernels}, 
are uniquely defined by the following properties. 
The Enriquez kernel $g^{i_{1} \cdots i_{r}}_{j}(x, y)$ is a $(1, 0)$-form in $x \in \tilde{R}$ 
and a scalar in $y \in \tilde{R}$ which is meromorphic for $x, y \in \tilde{R}$ 
and locally holomorphic in the complex moduli of pointed Riemann surfaces. 
In the above fundamental domain $D$, 
$g^{i_{1} \cdots i_{r}}_{j}(x, y)$ is holomorphic in $x, y \in D$ for $r \geq 2$, 
it can have a simple pole in $x$ at $y$ such as 
$g^{i}_{j}(x, y) - \frac{\delta_{ij} dx}{x - y}$ is holomorphic for $r = 1$, 
and is given by $g^{\emptyset}_{j}(x, y) = \omega_{j}(x)$ for $r = 0$. 
The monodromies in $x$ of $g^{i_{1} \cdots i_{r}}_{j}(x, y)$ around the $A$ cycles are trivial,
and those around the cycle $B_{l}$ are given by 
$$
g^{i_{1} \cdots i_{r}}_{j}(B_{l} \cdot x, y) = 
g^{i_{1} \cdots i_{r}}_{j}(x, y) + \sum_{k=1}^{r} \frac{(-2 \pi \sqrt{-1})^{k}}{k!} 
\delta^{i_{1} \cdots i_{k}}_{l} g^{i_{k+1} \cdots i_{r}}_{j}(x, y), 
$$
where $B_{l} \cdot x$ denotes the action of the element $B_{l} \in \pi_{1}(R, q)$ 
on the point $x \in R$, 
and the generalized Kronecker symbol is defined as 
$\delta^{i_{1} \cdots i_{k}}_{l} = \delta_{i_{1} l} \cdots \delta_{i_{k} l}$. 
Furthermore, by \cite{DESZ, DS}, 
these forms satisfy the following conditions. 

\begin{itemize} 

\item 
trivial $A$ monodromies in $y$, 
and $B$ monodromies given by 
$$
g^{i_{1} \cdots i_{r}}_{j}(x, B_{l} \cdot y) = 
g^{i_{1} \cdots i_{r}}_{j}(x, y) + 
\delta^{i_{r}}_{j} \sum_{k=1}^{r} \frac{(2 \pi \sqrt{-1})^{k}}{k!} 
g^{i_{1} \cdots i_{r-k}}_{j}(x, y) \delta^{i_{r-k+1} \cdots i_{r-1}}_{l}. 
$$

\item 
the periods around $A$ cycles on the boundary of the fundamental domain $D$ 
or any $y$ in the interior of $D$ are given in terms of Bernoulli numbers $B_{r}$ as 
$$
\oint_{A_{k}} g^{i_{1} \cdots i_{r}}_{j}(t, y) = 
(-2 \pi \sqrt{-1})^{r} \frac{B_{r}}{r!} \delta^{i_{1} \cdots i_{r} k}_{j}. 
$$

\end{itemize}
Then the ${\mathfrak g}$-valued Enriquez connection $d - {\mathcal K}_{R}$ is expressed as 
$$
{\mathcal K}_{R}(x, y) = \sum_{r = 0}^{\infty} \sum_{1 \leq i_{1},..., i_{r}, j \leq g} 
g^{i_{1} \cdots i_{r}}_{j} (x, y) \ [b_{i_{1}}, [b_{i_{2}}, \cdots [b_{i_{r}}, a_{j}] \cdots ]]. 
\eqno(2.1) 
$$ 

We recall the classical definition of the prime form following Bogatyr\"{e}v \cite{Bog}. 
Let $(\omega_{R})_{x,y}$ denote the unique meromorphic Abelian differential on $R$ 
whose poles are simple at $x, y$ with residues $1, -1$ respectively 
such that $(\omega_{R})_{x,y}$ is normalized for $\{ A_{i} \}$ in the sense that 
$\oint_{A_{i}} (\omega_{R})_{x,y} = 0$ for all $i = 1,..., g$. 
For points $\alpha, \beta$ on $R$ with local coordinates, 
take their lifts $\tilde{\alpha}, \tilde{\beta}$ on the universal cover $\tilde{R}$ of $R$ 
which give rise to a homotopy class of paths from $\beta$ to $\alpha$. 
Then the associated prime form is defined as 
$$
E (\tilde{\alpha}, \tilde{\beta}) = E_{R} (\tilde{\alpha}, \tilde{\beta}) = 
\frac{P (\tilde{\alpha}, \tilde{\beta})}{\sqrt{d \alpha} \sqrt{d \beta}}, 
$$
where 
$$
P (\tilde{\alpha}, \tilde{\beta}) = P_{R} (\tilde{\alpha}, \tilde{\beta}) := 
\left( \lim_{x \rightarrow \alpha, y \rightarrow \beta} 
\frac{-(\alpha - x)(\beta - y)}{\exp \left( \int_{\tilde{\beta}}^{\tilde{\alpha}} 
(\omega_{R})_{x, y} \right)} \right)^{1/2} 
$$
is continuous in $\alpha, \beta$ and satisfies 
$$
P (\tilde{\alpha}, \tilde{\beta}) = (\alpha - \beta)(1 + \mbox{higher order terms}) 
$$ 
for the local coordinate at $\alpha$ if $\beta$ is close to $\alpha$.
This implies 
$$
P (\tilde{\alpha}, \tilde{\beta}) = (\alpha - \beta) \exp 
\left( - \frac{1}{2} \int_{\tilde{\beta}}^{\tilde{\alpha}} (\omega_{R})^{*}_{x, y} \right); \ 
(\omega_{R})^{*}_{x, y} := (\omega_{R})_{x, y} - \left( \frac{dz}{z - x} - \frac{dz}{z - y} \right) 
\eqno(2.2)
$$
if $\alpha, \beta$ are close to each other.

Then \cite[Theorem 1]{DS} expresses the special kernel functions $g^{i}_{j}(x, y)$ 
in terms of Abelian differentials and prime forms: 
$$
g^{i}_{j}(x, y) = 
\oint_{A_{i}} \omega_{j}(t) \, \partial_{x} \log \left( \frac{E(x, y)}{E(x, t)} \right) + 
\pi \sqrt{-1} \delta^{i}_{j} \omega_{j}(x) 
\eqno(2.3) 
$$
with $x, y$ in the interior of $D$ and the cycle $A_{i}$ on the boundary of $D$. 
Furthermore, 
\cite[Theorem 2]{DS} gives recursion formulas 
between the general kernel functions  $g^{i_{1} \cdots i_{r}}_{j}(x, y)$: 
$$ 
\begin{array}{rcl}
g^{l i_{1} \cdots i_{r}}_{k}(y, z) 
& = & 
{\displaystyle - \oint_{A_{l}} \sum_{j=1}^{g} g^{j}_{k}(y, t) g^{i_{1} \cdots i_{r}}_{j}(t, z)} 
\\
& & 
{\displaystyle - \sum_{m=1}^{r-1} \sum_{l=1}^{g} (-2 \pi \sqrt{-1})^{m} \frac{B_{m}}{m!} 
\delta^{i_{1} \cdots i_{m}}_{l} g^{l i_{m+1} \cdots i_{r}}_{k}(y, z)} 
\\
& & 
{\displaystyle 
- \, \omega_{k}(y) (-2 \pi \sqrt{-1})^{r+1} \frac{B_{r+1}}{r!} \delta^{i_{1} \cdots i_{r} l}_{k}} 
\end{array}
\eqno(2.4)
$$
for $r \geq 1$ and $y, z$ in the interior of $D$ and cycles $A_{l}$ on the boundary of $D$. 
Furthermore, 
the convolution integral over $A_{l}$ in the first line of this formula is defined as a limit 
$$
\lim_{\varepsilon \rightarrow 0} \oint_{A_{l}^{\varepsilon}} 
\sum_{j=1}^{g} g^{j}_{k}(y, t) g^{i_{1} \cdots i_{r}}_{j}(t, z), 
$$
where the cycle $A_{l}^{\varepsilon}$ are small deformations of $A_{l}$ which are homotopic 
to $A_{l}$ and contained in the interior of $D$.

\section{Variation of the Enriquez connection} 

For each $i = 1,..., g$, 
let $(E_{i}; t_{i})$ be an elliptic (complex) curve, i.e., a Riemann surface of genus $1$ 
with one marked point $t_{i}$, 
and take a symplectic homology basis $\{ A_{i}, B_{i} \}$ of $E_{i} \setminus \{ t_{i} \}$. 
Denote by $R_{0}$ the singular complex curve obtained from 
the Riemann sphere $\mathbb{CP}^{1}$ with $g$ marked points $x_{1},..., x_{g}$ 
and the elliptic curves $(E_{1}; t_{1}),..., (E_{g}; t_{g})$ by identifying $x_{i}$ and $t_{i}$ $(i = 1,..., g)$. 
For complex parameters $s_{i}$ with small $|s_{i}|$, 
put $s = (s_{1},...,s_{g})$, 
and denote by $R_{s}$ the family of deformations of $R_{0}$ defined as 
$$
\xi_{i} \cdot \theta_{i} = s_{i} \quad (i = 1,..., g) 
\eqno(3.1)
$$
for local coordinates $\xi_{i}$ (resp. $\theta_{i}$) at $x_{i}$ (resp. $t_{i}$) such that 
$\xi_{i}(x_{i}) = 0$ (resp. $\theta_{i}(t_{i}) = 0$). 
Then one can take neighborhoods $U_{i}$ of $x_{i} \in \mathbb{CP}^{1}$ and 
$V_{i}$ of $t_{i} \in E_{i}$ $(i = 1,..., g)$ such that $R_{s}$ is obtained from $R_{0}$ by identifying 
the boundaries $\partial U_{i}, \partial V_{i}$ of $U_{i}, V_{i}$ respectively via (3.1). 
In other words, 
$R_{s}$ is obtained by gluing arounds of $\partial U_{i}$ and $\partial V_{i}$ via (3.1) such that 
the inside (resp. outside) of $\partial U_{i}$ are mapped to the outside (resp. inside) of $\partial V_{i}$. 
Therefore, if $s \rightarrow {\bf 0}:= (0,...,0)$, 
then (3.1) become the local equations $\xi_{i} \cdot \theta_{i} = 0$ 
at the singular points $x_{i} = t_{i}$ on $R_{0}$ which are called nodal (see {\bf Figure 1}). 
\begin{figure}
\centering 
\includegraphics[width=1.00\textwidth,origin=c]{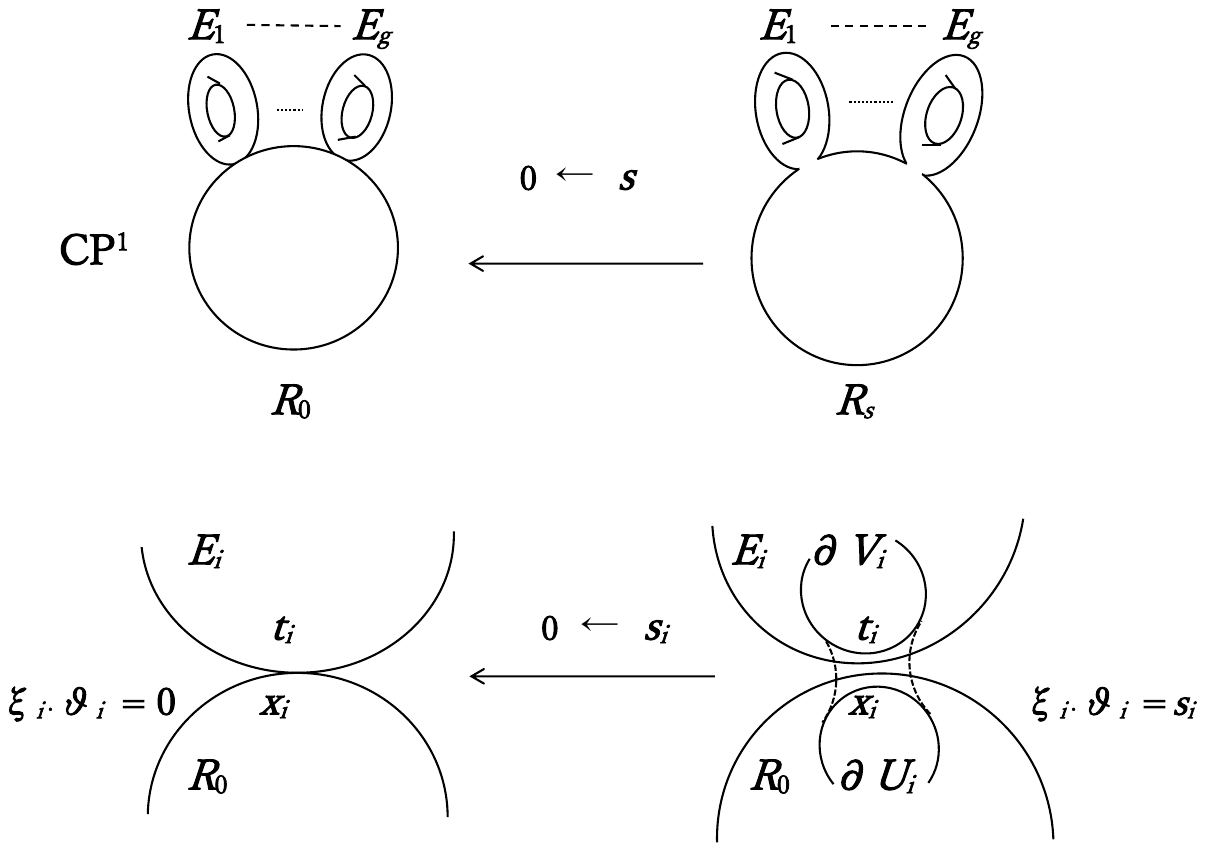}
\caption{\label{Figure 1} The singular curve $R_{0}$, 
its deformation $R_{s}$ and their local description.}
\end{figure}
The above $\{ A_{i}, B_{i} \}_{1 \leq i \leq g}$ 
gives a symplectic homology basis of $R_{s}$, 
and points on $E_{i} \setminus \{ t_{i} \}$ and on $\mathbb{CP}^{1} \setminus \{ x_{1},..., x_{g} \}$ 
are identified with those on $R_{s}$. 

First, 
we study the variation of Abelian differentials using the theory of stable Abelian differentials 
on semi-stable (algebraic) curves \cite[Section 1]{DM} 
(see also \cite[Chapter 3, Section A]{HM}). 
Similar arguments were given in \cite[Theorems 7.1 and 7.4]{I2}, 
and analytic approaches were given by Fay \cite{F} and Hu-Norton \cite{HN}. 
Recall that a projective (i.e., compact when the base field is $\mathbb{C}$) 
algebraic curve $C$ is called semi-stable when their singularities are nodal, 
and one can take the desingularization $C' \rightarrow C$ such that 
each singular point $p$ on $C$ are decomposed to $2$ smooth points $p_{1}, p_{2}$ on $C'$. 
Then a $1$-form on a semi-stable curve $C$ is called a stable Abelian differential 
when for each singular point $p$ on $C$, 
this form on $C'$ can have poles at $p_{1}, p_{2}$ which are simple such that 
the sum of the residues at $p_{1}, p_{2}$ is $0$. 
It is shown in \cite{DM, HM} that the spaces of stable Abelian differentials on semi-stable curves 
satisfy similar dimension formulas containing the Riemann-Roch theorem 
as for those of Abelian differentials on projective smooth curves 
(i.e., compact Riemann surfaces when the base field is $\mathbb{C}$). 

For each $i = 1,..., g$, 
there exists a unique family $(\omega_{s})_{i}$ of holomorphic Abelian differentials 
on the above $R_{s}$ with $s_{1},..., s_{g} \neq 0$ which is normalized for $\{ A_{i} \}$, 
i.e., $\oint_{A_{i}} (\omega_{s})_{j} = \delta_{ij}$. 
Furthermore, 
there exists a unique stable Abelian differential $(\omega_{0})_{i}$ on $R_{0}$ such that 
its restriction $(\omega_{0})_{i}|_{\mathbb{CP}^{1}}$ to $\mathbb{CP}^{1}$ is $0$, 
$(\omega_{0})_{i}|_{E_{j}} = 0$ for $j \neq i$ and $(\omega_{0})_{i}|_{E_{i}}$ 
is the unique holomorphic Abelian differential on $E_{i}$ normalized for $A_{i}$. 
It is shown in \cite{HM} that the sheaf of stable Abelian differentials on $R_{s}$ 
is locally free on $s$ of rank $g$ which states that there exists a set 
$\{ (\omega'_{s})_{i} \}_{1 \leq i \leq g}$ of stable Abelian differentials 
on small deformations $R_{s}$ of $R_{0}$ satisfying 
\begin{itemize}

\item 
each $(\omega'_{s})_{i}$ is holomorphic as a function of $s_{1},..., s_{g}$, 
and it tends to $(\omega_{0})_{i}$ under $s \rightarrow {\bf 0}$, 

\item 
the set $\{ (\omega'_{s})_{i} \}_{i}$ gives a basis of the space of stable Abelian differentials on $R_{s}$. 

\end{itemize}
Therefore, 
$$
\left( (\omega'_{s})_{1},..., (\omega'_{s})_{g} \right) \left( \oint_{A_{i}} (\omega'_{s})_{j} \right)_{i,j}^{-1} 
= \left( (\omega_{s})_{1},..., (\omega_{s})_{g} \right)
$$ 
and 
$$
\oint_{A_{i}} (\omega'_{s})_{j} \rightarrow \oint_{A_{i}} (\omega_{0})_{j}|_{E_{i}} = \delta_{ij} 
\quad \mbox{under $s \rightarrow {\bf 0}$} 
$$ 
which imply that 
$$
(\omega_{s})_{i} \rightarrow (\omega_{0})_{i} \quad \mbox{under $s \rightarrow {\bf 0}$.} 
\eqno(3.2)
$$ 
For points $x \in E_{i} \setminus \{ t_{i} \}$ with fixed $i$ and 
$y \in \mathbb{CP}^{1} \setminus \{ x_{1},..., x_{g} \}$, 
denote by $(\omega_{s})_{x, y}$ the unique family of meromorphic Abelian differentials on $R_{s}$ 
with $s_{1},..., s_{g} \neq 0$ whose poles are simple at $x, y$ with residues $1, -1$ respectively such that $(\omega_{s})_{x, y}$ is normalized, 
i.e., $\oint_{A_{i}} (\omega_{s})_{x, y} = 0$ for all $i = 1,..., g$. 
Furthermore, 
denote by $(\omega_{0})_{x, y}$ the unique stable Abelian differential on $R_{0}$ 
such that 
$(\omega_{0})_{x,y}|_{\mathbb{CP}^{1}} = \left( \frac{1}{z - x_{i}} - \frac{1}{z - y} \right) dz$, 
$(\omega_{0})_{x,y}|_{E_{j}} = 0$ for $j \neq i$ and $(\omega_{0})_{x,y}|_{E_{i}}$ 
is the unique meromorphic Abelian differential on $E_{i}$ normalized for $A_{i}$ 
whose poles are simple at $x, t_{i}$ with residues $1, -1$ respectively. 
If we represent $E_{i}$ as the quotient space of $\mathbb{C} \setminus \{ 0 \}$ 
by the action $z \mapsto q_{i} z$ with $0 < |q_{i}| < 1$, 
and take the cycle $A_{i}$ as a circle around $z = 0$, 
then applying general formulas of such Abelian differentials \cite[Section 3]{I1} 
to the genus $1$ case, 
one has 
\begin{eqnarray*}
(\omega_{0})_{x,y}|_{E_{i}}(z) 
& = & 
\sum_{n = -\infty}^{\infty} 
\left( \frac{1}{q_{i}^{n} z - x} - \frac{1}{q_{i}^{n} z - t_{i}} \right) d(q_{i}^{n} z) 
\\
& = & 
\left( \frac{1}{z - x} - \frac{1}{z - t_{i}} \right) dz 
\\
& & 
+ \, 
\sum_{n = 1}^{\infty} \left( \frac{q_{i}^{n}}{q_{i}^{n} z - x} - \frac{q_{i}^{n}}{q_{i}^{n} z - t_{i}} \right) dz 
+ 
\sum_{n = 1}^{\infty} 
\left( \frac{1}{z - q_{i}^{n} x} - \frac{1}{z - q_{i}^{n}t_{i}} \right) dz. 
\end{eqnarray*} 
Since the sheaf of meromorphic stable Abelian differentials on $R_{s}$ which have no pole 
or simple poles at $x, y$ is locally free on $s$ of rank $g+1$ (cf. \cite{HM}), 
there exists a family $(\omega'_{s})_{x, y}$ of meromorphic stable Abelian differentials 
on small deformations $R_{s}$ of $R_{0}$ which has simple poles at $x, y$ with residues $1, -1$ respectively 
such that $(\omega'_{s})_{x, y}$ is holomorphic as a function of $s_{1},..., s_{g}$, 
and it tends to $(\omega_{0})_{x, y}$ under $s \rightarrow {\bf 0}$. 
Therefore, 
$$
(\omega'_{s})_{x, y} - 
\sum_{i=1}^{g} \left( \oint_{A_{i}} (\omega'_{s})_{x, y} \right) (\omega_{s})_{i} = (\omega_{s})_{x, y}, 
$$
and hence 
$$
(\omega_{s})_{x, y} \rightarrow (\omega_{0})_{x, y} \quad \mbox{under $s \rightarrow {\bf 0}$.}
\eqno(3.3) 
$$
 
Second, 
we study the variation of the Enriquez kernels and connection under $s \rightarrow {\bf 0}$. 
Here we recall the explicit formula (2.3) of the Enriquez kernel $g_{j}^{i}$. 
If $i \neq j$ and $s \rightarrow {\bf 0}$, 
then (3.2) implies that $(\omega_{s})_{j}(t) \rightarrow 0$ for $t \in A_{i}$, 
and hence $(g_{R_{s}})_{j}^{i}(x, y) \rightarrow 0$. 
If $x \in E_{k} \setminus \{ t_{k} \}$ for $k \neq i$ and $s \rightarrow {\bf 0}$, 
then paths from $x$ to $y \in \mathbb{CP}^{1} \setminus \{ x_{1},..., x_{g} \}$ and to $t$ 
become those passing through $t_{k} = x_{k}$, 
and hence by (3.3) and (2.2), 
$$
\partial_{x} \log \left( \frac{E_{R_{s}} (x, y)}{E_{R_{s}} (x, t)} \right) \rightarrow 
\partial_{x} \log \left( \frac{E_{R_{0}} (x_{k}, y)}{E_{R_{0}} (x_{k}, t)} \right) = 0 
$$
which combined with $\delta_{j}^{i} (\omega_{s})_{j}(x) \rightarrow 0$ imply that 
$(g_{R_{s}})_{j}^{i}(x, y) \rightarrow 0$. 
Furthermore, if $x \in E_{i} \setminus \{ t_{i} \}$, 
then by (3.3), 
$$
\partial_{x} \log E_{R_{s}}(x, y) \rightarrow \partial_{x} \log E_{E_{i}}(x, t_{i}) \quad 
\mbox{under $s \rightarrow {\bf 0}$}, 
$$
and hence by (3.2), 
$$
(g_{R_{s}})_{j}^{i} (x, y) \rightarrow (g_{E_{i}})_{i}^{i}(x, t_{i}) \quad 
\mbox{under $s \rightarrow {\bf 0}$}. 
$$
Therefore, 
if $x \in E_{k} \setminus \{ t_{k} \}$ and $y \in \mathbb{CP}^{1} \setminus \{ x_{1},..., x_{g} \}$, 
then 
$$
(g_{R_{s}})_{j}^{i} (x, y) \rightarrow 
\left\{ \begin{array}{ll} 
(g_{E_{k}})_{k}^{k}(x, t_{k}) & (\mbox{if $i = j = k$}), 
\\ 
0 & (\mbox{otherwise})  
\end{array} \right. 
\quad \mbox{under $s \rightarrow {\bf 0}$}. 
\eqno(3.4)
$$
Take $x \in E_{i} \setminus \{ t_{i} \}$ and $y \in \mathbb{CP}^{1} \setminus \{ x_{1},..., x_{g} \}$. 
Then (3.4) and the recursion formulas (2.4) for $R_{s}$ and $E_{i}$ imply that 
$\lim_{s \rightarrow {\bf 0}} (g_{R_{s}})_{j}^{k i_{1}} (x, y)$ is equal to  
$$
- \oint_{A_{i}} (g_{E_{i}})_{i}^{i}(x, t) (g_{E_{i}})_{i}^{i}(t, t_{i}) 
- (\omega_{0})_{i}|_{E_{i}}(x) (-2 \pi \sqrt{-1})^{2} B_{2}  
= (g_{E_{i}})_{i}^{i i }(x, t_{i})
$$
if $i = j = k = i_{1}$, 
and it is $0$ otherwise. 
Similarly, by induction on $r$, 
$\lim_{s \rightarrow {\bf 0}} (g_{R_{s}})_{j}^{k i_{1} \cdots i_{r}} (x, y)$ is equal to 
\begin{eqnarray*}
- \oint_{A_{i}} (g_{E_{i}})_{i}^{i}(x, t) (g_{E_{i}})_{i}^{i_{1} \cdots i_{r}}(t, t_{i}) 
& - & \sum_{m=1}^{r-1} (-2 \pi \sqrt{-1})^{m} \frac{B_{m}}{m!} (g_{E_{i}})_{i}^{i i_{m+1} \cdots i_{r}}(x, y) 
\\ 
& - & (\omega_{0})_{i}|_{E_{i}}(x) (-2 \pi \sqrt{-1})^{r+1} \frac{B_{r+1}}{r!} 
\end{eqnarray*}
if $i = j = k = i_{1} = \cdots = i_{r}$, and it is $0$ otherwise, 
and hence we obtain the variational formula 
for $x \in E_{i} \setminus \{ t_{i} \}, y \in \mathbb{CP}^{1} \setminus \{ x_{1},..., x_{g} \}$:  
$$
\lim_{s \rightarrow {\bf 0}} (g_{R_{s}})_{j}^{k i_{1} \cdots i_{r}} (x, y) = 
\left\{ \begin{array}{ll} 
(g_{E_{i}})_{i}^{i i \cdots i}(x, t_{i}) & (\mbox{if $i = j = k = i_{1} = \cdots = i_{r}$}), 
\\ 
0 & (\mbox{otherwise}). 
\end{array} \right. 
$$
Therefore, 
if $x \in E_{i} \setminus \{ t_{i} \}, y \in \mathbb{CP}^{1} \setminus \{ x_{1},..., x_{g} \}$, 
then under $s \rightarrow {\bf 0}$, 
the connection form ${\mathcal K}_{R_{s}}(x, y)$ on $R_{s}$ given in (2.1) becomes ${\mathcal K}_{E_{i}} (x, t_{i})$ 
which is shown to be the elliptic KZB connection form 
which has a simple pole at $x = t_{i}$ with residue $[b_{i}, a_{i}]$ 
(cf. \cite[Section 8]{E1} and \cite[Section 4.2]{DS}). 
Furthermore, ${\mathcal K}_{R_{s}}(x, y)$ for 
$x, y \in \mathbb{CP}^{1} \setminus \{ x_{1},..., x_{g} \}$ becomes the KZ connection form 
${\mathcal K}_{\mathbb{CP}^{1}}$ on $\mathbb{CP}^{1}$ with trivial underlying bundle 
which has simple poles at $x = x_{i}$ with residues $-[b_{i}, a_{i}]$ $(i = 1,..., g)$, 
and $x = y$ with residue $\sum_{i=1}^{g} [b_{i}, a_{i}]$.

\section{Polylogarithm sheaves and functions} 

In the previous section, we showed that under $s = (s_{i}) \rightarrow {\bf 0}$, 
the Enriquez connection forms ${\mathcal K}_{R_{s}}$ on $R_{s}$ tend to 
the elliptic KZB connection forms ${\mathcal K}_{E_{i}}$ on $E_{i}$ $(i = 1,..., g)$ 
having simple poles at $t_{i}$ with residue $[b_{i}, a_{i}]$, 
and to the KZ connection form ${\mathcal K}_{\mathbb{CP}^{1}}$ (with trivial bundle) 
on  $\mathbb{CP}^1$ having simple poles at $x_{i}$ with residue $-[b_{i}, a_{i}]$ and 
at $y$ with residue $\sum_{i=1}^{g} [b_{i}, a_{i}]$. 
Using this result, 
we will describe the Enriquez connection on $R_{s}$ as the {\it polylogarithm sheaf} given in 
\cite[Section 3.2]{I3}) for $g \geq 2$. 

First, 
we consider a degeneration of this $\mathbb{CP}^{1}$ 
with $g+1$ marked points $x_{1},..., x_{g}, y$ to a union $P_{0}$ of 
$g-1$ copies of $\mathbb{CP}^{1}$ with $3$ normalized marked points $0, 1, \infty$. 
Then the dual graph $\Delta$ of $P_{0}$ is a trivalent tree which consists of 
the sets $V$ of $g-1$ vertices, $E$ of $g-2$ oriented edges and $T$ of $g+1$ tails 
corresponding to the irreducible components $P_{v} = \mathbb{CP}^{1}$ $(v \in V)$ 
of $P_{0}$, 
the singular points connecting them and the marked points $x_{1},..., x_{g}, y$ respectively. 
Put $\pm E = \{ \pm e \mid e \in E \}$, 
where $-e$ denotes the inverse edge of $e$, 
and for each $h \in \pm E \cup T$, 
denote by $v(h)$ the target (resp. boundary) vertex of $h$ if $h \in \pm E$ (resp. $\in T$), 
and by $p_{h}$ the singular (resp. marked) point of $P_{v(h)}$ associated with $h$ 
if $h \in \pm E$ (resp. $\in T$). 
Therefore, 
for a given trivalent tree $\Delta$, 
the associated singular curve $P_{0}$ as a semi-stable curve obtained from 
$P_{v} = \mathbb{CP}^{1}$ $(v \in V)$ with $3$ normalized marked points $0, 1, \infty$ 
by identifying the points $p_{e} \in P_{v(e)}$ and $p_{-e} \in P_{v(-e)}$ $(e \in E)$. 
In the case that $g = 3$, 
we give an example of $\Delta$ and $P_{0}$ (see {\bf Figure 2}): 
\begin{eqnarray*}
\Delta & = & (V = \{ v_{1}, v_{2} \}, \, E = \{ e \}, \, T = \{ t_{1}, t_{2}, t_{3}, t_{4} \}), 
\\
P_{0} & = & P_{v_{1}} \bigcup_{p_{e} = p_{-e}} P_{v_{2}}; \, 
\left\{ \begin{array}{l} 
P_{v_{1}} = \mathbb{CP}^{1} \ni p_{t_{1}} = 1, \, p_{t_{2}} = \infty, \, p_{-e} = 0, 
\\
P_{v_{2}} = \mathbb{CP}^{1} \ni p_{t_{3}} = 1, \, p_{t_{4}} = \infty, \, p_{e} = 0. 
\end{array} \right. 
\end{eqnarray*}

\begin{figure}[tbp]
\centering 
\includegraphics[width=0.950\textwidth,origin=c]{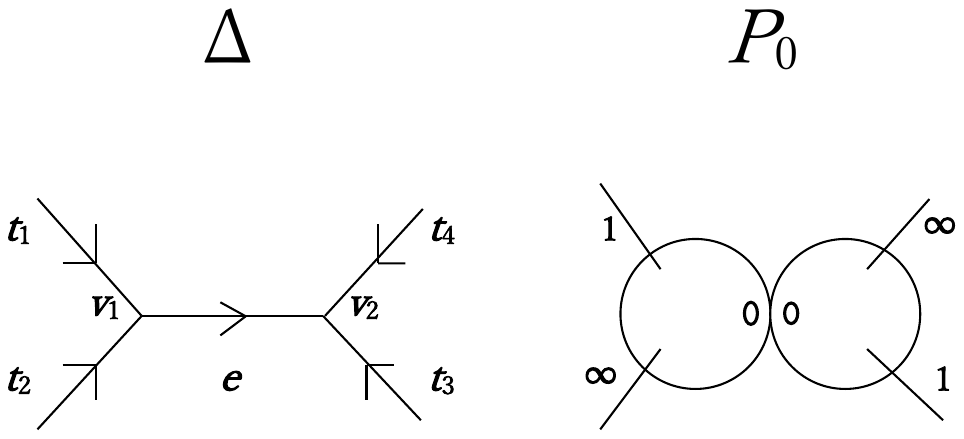}
\caption{\label{Figure 2} 
A trivalent tree $\Delta$ and the associated singular curve $P_{0}$.}
\end{figure}

\noindent
Under this degeneration $\mathbb{CP}^{1} \rightarrow P_{0}$, 
the above KZ connection form ${\mathcal K}_{\mathbb{CP}^{1}}$ is reduced to 
the basic KZ connection forms on $P_{v}$ $(v \in V)$ denoted by 
$\left( \frac{R_{v,0}}{z} + \frac{R_{v,1}}{z-1} \right) dz$ which have simple poles at $0, 1$ and $\infty$ 
with residues $R_{v,0}, R_{v,1}$ and $R_{v,\infty} = - R_{v,0} - R_{v,1}$ respectively. 
Furthermore, 
the sum of $2$ residues at the desingularization of each singular point on $P_{0}$ is $0$. 
Therefore, 
the residues $X_{h}$ at $p_{h} \in P_{v(h)}$ $(h \in \pm E \cup T)$ 
of the basic KZ connections satisfy the relations: 
\begin{itemize}

\item 
if $t \in T$ corresponds to $x_{i}$ (resp. $y$), 
then $X_{t}$ is $-[b_{i}, a_{i}]$ (resp. $\sum_{i=1}^{g} [b_{i}, a_{i}]$), 

\item 
if $h$ is an oriented edge $h$ of $\Delta$, then $X_{h} + X_{-h}$ is $0$, 

\item 
if $v \in V$, 
then the sum of $X_{h}$ for oriented edges $h$ with terminal vertex $v$ is $0$, 

\end{itemize}
and hence each $X_{h}$ is uniquely represented as a sum of $[b_{i}, a_{i}]$ $(i = 1,..., g)$ 
with integral coefficients. 
In the above example when $g = 3$, 
if $x_{i}$ correspond to $p_{t_{i}}$ $(i = 1, 2, 3)$ and $y$ corresponds to $p_{t_{4}}$, 
then 
$$
X_{t_{i}} = - [b_{i}, a_{i}] \ (i = 1, 2, 3), \quad X_{t_{4}} = \sum_{i=1}^{3} [b_{i}, a_{i}], 
$$
and hence 
$$
X_{-e} = [b_{1}, a_{1}] + [b_{2}, a_{2}], \quad X_{e} = - [b_{1}, a_{1}] - [b_{2}, a_{2}]. 
$$
Following \cite[Section 3]{I3}, 
we recover the KZ connection on small deformations of $P_{0}$ by gluing the basic KZ connections 
on $P_{v}$ $(v \in V)$ with residues $X_{h}$ for $h \in \pm E$ satisfying $v(h) = v$. 
For each marked point $p_{h}$ corresponding to $h \in \pm E$, 
fix its local coordinate $z_{h}$ given by 
$$
z_{h} = \left\{ \begin{array}{ll} 
z - p_{h} & (\mbox{if $p_{h} = 0, 1$}), 
\\ 
1/z & (\mbox{if $p_{h} = \infty$}), 
\end{array} \right. 
$$
and define the deformation of $P_{0}$ by parameters $s_{h} = s_{-h}$ $(h \in \pm E)$ via the relations 
$$
z_{h} \cdot z_{-h} = s_{h} \quad (h \in \pm E) \eqno(4.1)
$$ 
which is identified with $\mathbb{CP}^{1}$ obtained by patching $P_{v(h)} \setminus U_{h}$ to 
$U_{-h}$ $(h \in \pm E)$ for certain open neighborhoods $U_{h}$ of $p_{h} \in P_{v(h)}$. 
Since $X_{h} dz_{h}/z_{h} = X_{-h} dz_{-h}/z_{-h}$, 
one can glue the basic KZ connections on $P_{v(h)}$ and $P_{v(-h)}$ $(h \in \pm E)$, 
and obtain the KZ connection on $\mathbb{CP}^{1}$ with simple poles at $x_{i}$ with residue $-[b_{i}, a_{i}]$ 
and at $y$ with residue $\sum_{i=1}^{g} [b_{i}, a_{i}]$.

Second, 
we show that the Enriquez connection becomes the polylogarithm sheaf obtained by gluing 
the basic KZ connections and the elliptic KZB connections. 
We take the above representation of the elliptic curve $(E_{i}; t_{i})$ as 
$\left( (\mathbb{C} \setminus \{ 0 \})/\langle q_{i} \rangle; 1 \right)$ with $0 < |q_{i}| < 1$, 
and construct $R_{s}$ by gluing $P_{v}$ $(v \in V)$ and $E_{i}$ $(i = 1,..., g)$ by the relations (4.1) and (3.1), 
where $\xi_{i} = z_{h_{i}}$ if $p_{h_{i}} = x_{i}$ and $\theta_{i} = z - 1$. 
Then by the equation $-[b_{i}, a_{i}] d \xi_{i}/\xi_{i} = [b_{i}, a_{i}] d \theta_{i}/\theta_{i}$,  
the basic KZ connections on $P_{v}$ and the elliptic KZB connections on $E_{i}$ can be glued to 
a connection on $R_{s}$ which is the Enriquez connection. 
By this description, 
one can see that the polylogarithm sheaf originally defined on the local family $R_{s}$ 
of pointed Riemann surfaces can be extended to the Enriquez connection on these whole family. 
Furthermore, 
monodromies of the Enriquez connection were expressed explicitly in \cite[Section 3.3]{I3} 
as noncommutative formal power series in $X_{h}$ $(h \in \pm E \cup T)$ whose coefficients are 
power series in deformation parameters and their logarithms with coefficients given by multiple zeta values. 
Therefore, we obtain explicit formulas of the higher genus polylogarithms for the Enriquez connection. 
Since the Enriquez connection is defined for any Riemann surfaces, 
these explicit formulas can be analytically continued for these global families whose calculation 
was also given using the Schottky uniformization theory \cite[Section 3.3 and Appendix A]{I3}

\section{Conclusion and outlook}
We studied the variation of the Enriquez connections \cite{E1} 
under degenerations of pointed Riemann surfaces, 
and showed that the connections can be described as the polylogarithm sheaves given in \cite{I3}. 
Therefore, 
one can see that the polylogarithm sheaves defined on local families of pointed Riemann surfaces 
can be globally extended to the whole families, 
and that the higher genus polylogarithms derived from the Enriquez connection 
can be calculated explicitly. 
We hope that our explicit formulas of the polylogarithms will be applicable to expressing hyperlogarithms 
on families of pointed Riemann surfaces which were studied in the genus $1$ case by Br\"{o}del 
and others \cite{BDDT, BDDPT1, BDDPT2, BDDPT3, BK} in order to calculate Feynman integrals. 

Note that the results of Section 3 for the variation of Abelian differentials on Riemann surfaces 
can be extended to general degenerations which contain separating and nonseparating ones 
(see \cite[Theorems 7.1 and 7.4]{I2}. 
Therefore, 
it is hoped that the variation of the Enriquez connections can be described in terms of 
the {\it tropical data} associated with considering degeneration as is done in 
\cite[Theorem 7.2]{I2} for the period integrals $\Omega_{ij} = \oint_{B_{i}} \omega_{j}$ 
of the Abelian differentials $\omega_{j}$. 
Furthermore, 
by the Schottky uniformization theory on families of Riemann surfaces 
with general degeneration \cite[Section 2]{I1}, 
we may obtain more precise variational formulas of the Enriquez connection 
combining explicit formulas of Abelian differentials in \cite[Section 3]{I1}
with the results on the Enriquez kernels in \cite{DS} referred above and 
in Baune-Broedel-Im-Lisitsyn-Zerbini \cite{BBILZ}. 
\vspace{4ex} 

\noindent 
{\bf Declaration of competing interest} 
\vspace{1ex}

The author declares that he has no known competing financial interests 
or personal relationships that could have appeared to influence the work 
reported in this article. 
\vspace{2ex}

\noindent 
{\bf Data availability} 
\vspace{1ex}

No data was used for the research described in this article. 
\vspace{2ex}

\noindent 
{\bf Acknowledgments} 
\vspace{1ex}

The author would like to express deep thanks to the referee 
whose comments were helpful in modifying this manuscript.  
This work is partially supported by the JSPS Grant-in-Aid for 
Scientific Research No. 25K06920.





\end{document}